\documentclass{article}%
\usepackage{utopia}
\usepackage{amssymb,latexsym}
\usepackage{amsfonts}
\usepackage{amsmath}
\usepackage[colorlinks=true, pdfstartview=FitV, linkcolor=blue, citecolor=red,
urlcolor=blue]{hyperref}
\usepackage{color}
\usepackage{amssymb}
\usepackage{graphicx}%
\providecommand{\U}[1]{\protect\rule{.1in}{.1in}}
\newtheorem{theorem}{Theorem}

\newtheorem{corollary}[theorem]{Corollary}

\allowdisplaybreaks
\begin{document}

\title{Formulas for sums of powers of integers and their reciprocals}
\author{Levent Karg\i n\thanks{leventkargin48@gmail.com}, Ayhan
Dil\thanks{adil@akdeniz.edu.tr}, M\"{u}m\"{u}n Can\thanks{mcan@akdeniz.edu.tr}\\Department of Mathematics, Akdeniz University, Antalya, Turkey}
\date{}
\maketitle

\begin{abstract}
This paper gives new explicit formulas for sums of powers of integers and
their reciprocals.



\end{abstract}

\markright{Formulas for sums of powers}

\section{Introduction}

For an integer $p>0$, consider the following sum of powers of integers
\begin{equation}
\sum_{k=1}^{n}k^{p}=1^{p}+2^{p}+\cdots+n^{p}. \label{ps}%
\end{equation}
These sums have been of interest to mathematicians since antiquity. Over the
years, mathematicians have given formulas for special values of $p$. Johann
Faulhaber (1580-1635) proved that the consecutive powers of integers can be
expressed as a polynomial in $n$ of degrees $\left(  p+1\right)  $ and gave a
calculation formula up to $p=17$. The general form was established with the
discovery of the Bernoulli numbers $B_{n}$ as:%
\[
\sum_{k=1}^{n}k^{p}=\frac{1}{p+1}\sum_{k=0}^{p}\left(  -1\right)  ^{k}%
\binom{p+1}{k}B_{k}n^{p+1-k}.
\]
After that this formula is named as Faulhaber's formula. More details about
this formula can be found in \cite{Knuth} and the references therein. In 1978,
Gould \cite{G} gave an explicit formula for this sum as
\[
\sum_{k=1}^{n}k^{p}=\sum_{j=0}^{p}j!%
\genfrac{\{}{\}}{0pt}{}{p}{j}%
\binom{n+1}{j+1}=\sum_{j=0}^{p}\left(  -1\right)  ^{p+j}j!%
\genfrac{\{}{\}}{0pt}{}{p}{j}%
\binom{n+j}{j}.
\]
Recently, Merca \cite{Merca} expressed this sum in terms of the Stirling
numbers of the first and second kind. In this manner from past to present,
studies are ongoing to give new calculation formulas and also new proofs for
known formulas of this type of sums.

Besides, the sum of reciprocals of powers of the first $n$ natural numbers
corresponds to
\begin{equation}
\sum_{k=1}^{n}\frac{1}{k^{p}}=1+\frac{1}{2^{p}}+\frac{1}{3^{p}}+\cdots
+\frac{1}{n^{p}},\text{ }p>0 \label{psr}%
\end{equation}
which is also interesting and has a long history dating back to Leonard Euler.

In another context, the number formed by the sum of (\ref{ps}) and (\ref{psr})
together is called the $n$th generalized harmonic number and denoted by
$H_{n}^{\left(  p\right)  }$, namely%
\begin{equation}
H_{n}^{\left(  p\right)  }=\sum_{k=1}^{n}k^{-p},\text{ }p\in\mathbb{Z}.
\label{ghn}%
\end{equation}
As can be easily understood from (\ref{ghn}), the case $p<0$ reduces to
(\ref{ps}). When $p>1$ the number $H_{n}^{\left(  p\right)  }$ is the $n$th
partial sum of the famous Riemann zeta function. Its close relationship with
Riemann zeta function makes generalized harmonic numbers valuable in the field
of analytical number theory.

The case $p=1,$ $H_{n}^{\left(  1\right)  }=H_{n}$, is known as harmonic
number which occurs in fundamental equations in many areas from analysis to
discrete mathematics and computer science \cite{GKP,CG}. Harmonic numbers have
various relationships with both Bernoulli and Stirling numbers, we would like
to remind you one of them \cite[p. 424]{CSE}:
\begin{equation}
H_{n+1}=\frac{1}{n!}\sum_{k=0}^{n}\left(  -1\right)  ^{k}%
\genfrac{[}{]}{0pt}{}{n+1}{k+1}%
B_{k}, \label{2}%
\end{equation}
where $%
\genfrac{[}{]}{0pt}{}{n}{k}%
$ denotes the Stirling numbers of the first kind, the number of permutations
of $n$ elements with $k$ disjoint cycles.

In this paper, two formulas are given for the sum of powers of positive
integers, and a formula for their reciprocals.

In the first result, we present a formula for the generalized harmonic
numbers. Special cases of which correspond to the sums (\ref{ps}) and
(\ref{psr}). Apart from this, it is a general form of (\ref{2}) and answers
the question of which type of Bernoulli numbers are related to the generalized
harmonic numbers.

\begin{theorem}
\label{teo1}For all integers $p$ and non-negative integers $n$, we have
\begin{equation}
H_{n+1}^{\left(  p\right)  }=\sum_{k=1}^{n+1}\frac{1}{k^{p}}=\frac{1}{n!}
\sum_{j=0}^{n}
\genfrac{[}{]}{0pt}{}{n+1}{j+1}
B_{j}^{\left(  p\right)  }. \label{hpb}%
\end{equation}

\end{theorem}

Here, $B_{k}^{\left(  p\right)  }$ is the $k$th\ poly-Bernoulli number defined
by \cite[Eq. (1)]{Kaneko}
\begin{equation}
\sum_{k=1}^{\infty}B_{k}^{\left(  p\right)  }\frac{t^{k}}{k!}=\frac
{Li_{p}\left(  1-e^{-t}\right)  }{1-e^{-t}}, \label{pb}%
\end{equation}
where $Li_{p}\left(  z\right)  $ is the polylogarithm and has the generating
function
\[
Li_{p}\left(  z\right)  =\sum_{k=1}^{\infty}\frac{z^{k}}{k^{p}}.
\]
The poly-Bernoulli numbers are a generalization of the classical Bernoulli
number with $B_{k}^{\left(  1\right)  }=\left(  -1\right)  ^{k}B_{k}.$ They
have interesting combinatorial interpretations, and also appear in special
values of certain zeta functions.

Now, we present formulas for the sum of powers of integers, first of which is
as a consequence of Theorem \ref{teo1}.

\begin{corollary}
\label{cor}For all non-negative integers $n$ and $p$, we have
\[
\sum_{k=1}^{n}k^{p}=\sum_{j=0}^{p}j!
\genfrac{\{}{\}}{0pt}{}{p+1}{j+1}
\binom{n}{j+1}.
\]

\end{corollary}

\begin{theorem}
\label{teo2}For all positive integers $n$ and $p,$ we have%
\[
\sum_{k=1}^{n}k^{p}=\sum_{j=0}^{p}\left(  -1\right)  ^{p+j}j!
\genfrac{\{}{\}}{0pt}{}{p+1}{j+1}
\binom{n+j+1}{j+1}.
\]
\bigskip
\end{theorem}

Here, $%
\genfrac{\{}{\}}{0pt}{}{n}{k}%
$ is the Stirling numbers of the second kind, count the number of ways to
partition a set of $n$ objects into $k$ non-empty subsets.

\section{Proofs}

\subsection{Proof of Theorem \ref{teo1}}

Substituting $t\rightarrow1-e^{-t\text{ }}$ in the generating function of the
generalized harmonic numbers%
\begin{equation}
\sum_{k=0}^{\infty}H_{k}^{\left(  p\right)  }t^{k}=\frac{Li_{p}\left(
t\right)  }{1-t},\text{ }\left\vert t\right\vert <1, \label{gfgh}%
\end{equation}
one can easily obtain that%
\[
\frac{Li_{p}\left(  1-e^{-t}\right)  }{1-e^{-t}}=\sum_{n=0}^{\infty}\left(
-1\right)  ^{n}H_{n+1}^{\left(  p\right)  }\left(  e^{-t}-1\right)  ^{n}%
e^{-t}.
\]
We now utilize the following generating function of the second kind Stirling
numbers \cite[p. 351]{GKP}
\[
\sum_{k=n}^{\infty}
\genfrac{\{}{\}}{0pt}{}{k+1}{n+1}
\frac{z^{k}}{k!}=\frac{\left(  e^{z}-1\right)  ^{n}}{n!}e^{z},
\]
and deduce that%
\[
\frac{Li_{p}\left(  1-e^{-t}\right)  }{1-e^{-t}}=\sum_{k=0}^{\infty}\left(
\sum_{n=0}^{k}\left(  -1\right)  ^{k-n}
\genfrac{\{}{\}}{0pt}{}{k+1}{n+1}
n!H_{n+1}^{\left(  p\right)  }\right)  \frac{t^{k}}{k!}.
\]
Considering (\ref{pb}) and equating coefficients of $t^{n}$ give
\[
B_{n}^{\left(  p\right)  }=\sum_{k=0}^{n}\left(  -1\right)  ^{n-k}
\genfrac{\{}{\}}{0pt}{}{n+1}{k+1}
k!H_{k+1}^{\left(  p\right)  }.
\]
Finally, taking $a_{n}=B_{n}^{\left(  p\right)  }$ and $b_{k}=k!H_{k+1}%
^{\left(  p\right)  }$ in the well-known Stirling transform \cite[p.
310]{GKP}
\[
a_{n}=\sum_{k=0}^{n}\left(  -1\right)  ^{n-k}
\genfrac{\{}{\}}{0pt}{}{n+1}{k+1}
b_{k}\text{ if and only if }b_{n}=\sum_{k=0}^{n}
\genfrac{[}{]}{0pt}{}{n+1}{k+1}
a_{k}
\]
give the desired result.

\subsection{Proof of Corollary \ref{cor}}

Since $H_{n}^{\left(  -p\right)  }=\sum_{k=1}^{n}k^{p},$ take $-p$ ($p>0$) in
(\ref{hpb}) and then utilize the following identity of poly-Bernoulli numbers
\cite[Theorem 2]{AK}%
\[
B_{k}^{\left(  -p\right)  }=\sum_{j=0}^{\min\left\{  k,p\right\}  }\left(
j!\right)  ^{2}
\genfrac{\{}{\}}{0pt}{}{p+1}{j+1}
\genfrac{\{}{\}}{0pt}{}{k+1}{j+1}
.
\]
We find that
\[
\sum_{k=1}^{n+1}k^{p}=\frac{1}{n!}\sum_{j=0}^{n}\left(  j!\right)  ^{2}
\genfrac{\{}{\}}{0pt}{}{p+1}{j+1}
\sum_{k=j}^{n}
\genfrac{[}{]}{0pt}{}{n+1}{k+1}
\genfrac{\{}{\}}{0pt}{}{k+1}{j+1}
.
\]
Hence, thanks to the identity \cite[Theorems 3.7 and 3.11]{NyRa}
\[
\sum_{k=j}^{n}
\genfrac{[}{]}{0pt}{}{n+1}{k+1}
\genfrac{\{}{\}}{0pt}{}{k+1}{j+1}
=\frac{n!}{j!}\binom{n+1}{j+1},
\]
the proof is completed.

\subsection{Proof of Theorem \ref{teo2}}

With the aid of $H_{n}^{\left(  -p\right)  }=\sum_{k=1}^{n}k^{p},$ $p>0,$ we
may consider (\ref{gfgh}) as
\[
\sum_{n=0}^{\infty}\left(  \sum_{k=1}^{n}k^{p}\right)  t^{n}=\frac
{Li_{-p}\left(  t\right)  }{\left(  1-t\right)  }
\]
and use \cite{W}%
\[
Li_{-p}\left(  t\right)  =\left(  -1\right)  ^{p+1}\sum_{k=0}^{p}k!
\genfrac{\{}{\}}{0pt}{}{p+1}{k+1}
\left(  \frac{-1}{1-t}\right)  ^{k+1},\text{ }p\geq1
\]
to conclude that%
\[
\sum_{n=0}^{\infty}\left(  \sum_{k=1}^{n}k^{p}\right)  t^{n}=\sum
_{n=0}^{\infty}\left(  \sum_{k=0}^{p}\left(  -1\right)  ^{p+k}k!
\genfrac{\{}{\}}{0pt}{}{p+1}{k+1}
\binom{k+1+n}{k+1}\right)  t^{n}
\]
which is what we wanted to prove.


\begin{thebibliography}{99}                                                                                               %


\bibitem {AK}Arakawa T., Kaneko M. (1999). On poly-Bernoulli numbers.
\textit{Comment. Math. Univ. St. Pauli} 48(2): 159--167.

\bibitem {CSE}Cheon G.-S., El-Mikkawy M. E. A. (2008). Generalized harmonic
numbers with Riordan arrays. \textit{J. Number Theory} 128(2): 413--425.

\bibitem {G}Gould H. W. (1978). Evaluation of sums of convolved powers using
Stirling and Eulerian numbers. \textit{Fibonacci Quart.} 16: 488--497.

\bibitem {GKP}Graham R. L., Knuth D. E., Patashnik O. (1994). \textit{Concrete
Mathematics}. New York: Addison-Wesley.

\bibitem {Kaneko}Kaneko M. (1997). Poly-Bernoulli numbers. \textit{J. Theor.
Nr. Bordx.} 9: 221--228.

\bibitem {Knuth}Knuth D. E. (1993). Johann Faulhaber and sums of powers.
\textit{Math. Comput.} 61(203): 277--294.

\bibitem {CG}Conway J. H., Guy R. K. (1996). \textit{The Book of Numbers}. New
York: Springer-Verlag.

\bibitem {Merca}Merca M. (2015). An alternative to Faulhaber's formula.
\textit{Amer. Math. Monthly.} 122(6): 599-601.

\bibitem {NyRa}Nyul G., R\'{a}cz G. (2015). The $r$-Lah numbers.
\textit{Discrete Math.} 338(10): 1660--1666

\bibitem {W}Wood D. C. (1992). The computation of polylogarithms. Technical
Report 15-92. Canterbury, UK: University of Kent Computing Laboratory.
\end{thebibliography}
\end{document}